\documentclass[12pt]{article}
\usepackage{mathrsfs}
\usepackage{amssymb} \textwidth 150mm \textheight 220mm \oddsidemargin
15pt \evensidemargin 0pt \topmargin 0cm \headsep 0.3cm

\usepackage{amsmath}
\usepackage{graphicx}
\usepackage{indentfirst}
\usepackage{cite}
\usepackage{float}
\usepackage{dsfont}
\usepackage{CJK}
\usepackage{multirow}
\usepackage{makecell}

\usepackage{subfigure}

\begin{document}
\title{\Large\bf{Bifurcation of limit cycles from a cubic reversible  isochrone \thanks{E-mail address: jihua1113@163.com} }}
\author{{Jihua Yang,\ Qipeng Zhang}
\\ {\small \it School of Mathematical Sciences, Tianjin Normal University, }\\
{\small \it   Tianjin 300387, People's Republic of China}\\
 }
\date{}
\maketitle \baselineskip=0.9\normalbaselineskip \vspace{-3pt}
\noindent
{\bf Abstract}\, This paper is devoted to study the limit cycle problem of a cubic reversible system with an isochronous center, when it is perturbed inside a class of polynomials. An upper bound of the number of limit cycles is obtained using the Abelian integral. The algebraic structure of the Abelian integral is acquired thanks to some iterative formulas, which differs in many aspects from other methods. Some numerical simulations verify the existence of limit cycles.
\vskip 0.2 true cm
\noindent
 \section{Introduction and main result}
 \setcounter{equation}{0}
\renewcommand\theequation{1.\arabic{equation}}

Consider a planar vector field
\begin{eqnarray}
\frac{dx}{dt}=\dot{x}=P(x,y),\quad \frac{dy}{dt}=\dot{y}=Q(x,y),
\end{eqnarray}
where $P(x,y)$ and $Q(x,y)$ are analytic functions. $O$ is a singular point of system (1.1) which we may place at the origin. We say that $O$ is a center if there is a neighborhood  of it where all the orbits are periodic.  If the periodic orbits within a neighborhood of a center have the same period, then the center is called an isochronous center and the system (1.1) is called an isochronous system. In this case, the center is referred to
as uniformly isochronous, if the angular velocity remains constant to complete a full revolution around the origin \cite{C,CR,Y,LI}. Isochronous systems are not rare, see \cite{CL06,CL07,CL08}. The research on isochronous centers began with Galileo's study of the pendulum problem. Isochronism appears in many practical problems, and it manifests as a synchronous oscillation phenomenon. In fact, the most striking examples of isochronicity in the natural realm appear to be biology, encompassing crucial phenomena of human life such as the regular beating of our hearts and the circadian clocks that reflect our biological tempi \cite{S}.

It is well known that the problem of isochronicity appears only for the non-degenerate centers. If system (1.1) is a polynomial system of degree $n$ with a non-degenerate center at the origin, then it takes the form
\begin{eqnarray}
\dot{x}=-y+\sum\limits_{i+j=2}^na_{i,j}x^iy^j,\ \dot{y}=x+\sum\limits_{i+j=2}^nb_{i,j}x^iy^j,\ a_{i,j},b_{i,j}\in\mathds{R},
\end{eqnarray}
in an appropriate coordinate system and upon a time rescaling. When system (1.2) is a linear plus quadratic or cubic homogeneous polynomial differential system, the classification of isochronous centers has been solved in \cite{L,P}. There are a lot of articles studying their limit cycles \cite{LLL,CLYZ,FG,CJ,GV,ELV,CNT,GLLZ}, period functions \cite{MV,GY,LLLW} and so on. % phase portraits \cite{LV,GCZ}.
For the linear plus cubic homogeneous polynomial differential system, Pleshkan \cite{P} attested that the origin is an isochronous center if and only if this system can be brought to one of the following systems:\begin{eqnarray*}
\begin{aligned}
&S^*_1:\ \dot{x}=-y+x^3-3xy^2,\ \ \ \dot{y}=x+3x^2y-y^3,\\
&S^*_2:\ \dot{x}=-y+x^3-xy^2,\ \ \ \dot{y}=x+x^2y-y^3,\\
&S^*_3:\ \dot{x}=-y+3x^2y,\ \ \ \dot{y}=x-2x^3+9xy^2,\\
&S^*_4:\ \dot{x}=-y-3x^2y,\ \ \ \dot{y}=x+2x^3-9xy^2.\\
\end{aligned}
\end{eqnarray*}
 Grau and Villadelprat \cite{GV} proved that at most two critical periods bifurcate from each period annulus of  $S^*_1$, $S^*_2$, $S^*_3$ and $S^*_4$ when they are perturbed by polynomials of degree 3, and this bound can be reached. The limit cycle problem of $S^*_1$ under perturbations of polynomials of degree $n$ was considered by Gasull et al. in \cite{GLLZ}. The authors acquire that this perturbed system has at most $3n-1$ limit cycles for $n\geq9$. We haven't seen more results on the bifurcation of limit cycles about the above four isochronous systems yet.

%Now we give isochronous centers of cubic systems reversible with respect to a line passing through the origin.
For a linear plus quadratic and cubic homogeneous polynomial differential system of the form
\begin{eqnarray}
\dot{x}=-y+\sum\limits_{i+j=2}^3a_{i,j}x^iy^j,\quad \dot{y}=x+\sum\limits_{i+j=2}^3b_{i,j}x^iy^j,
\end{eqnarray}
the relevant research literatures are limited. In polar coordinates $(x,y)=(r\cos\theta,\ r\sin\theta)$, the system (1.3) can be written as
\begin{eqnarray}
\begin{cases}
\dot{r}=&r^2\big(\alpha_3\sin3(\theta+\theta_0)+\alpha_1\sin(\theta+\theta_0)\big)\\&+\beta_4\big(\alpha_4\sin4(\theta+\theta_0)+\alpha_2\sin2(\theta+\theta_0)\big),\\
\dot{\varphi}=&1+r\big(\alpha_3\cos3(\theta+\theta_0)+\beta_2\cos(\theta+\theta_0)\big)\\&+r^2\big(\alpha_4\cos4(\theta+\theta_0)+\beta_3\cos2(\theta+\theta_0)+\beta_1\big),\\
\end{cases}
\end{eqnarray}
here $\theta_0$, $\alpha_i$ and $\beta_i$ $(i=1,2,3,4)$ can be expressed by $a_{i,j}$ and $b_{i,j}$. %  and can be chosen as free parameters.
A general theory regarding isochronous problem of system (1.4) has not yet been developed. There is a conclusion when system (1.4) is a reversible system (with respect to a straight line passing through the origin) for $\alpha_3=\alpha_4=0$, which can be found  in \cite{CG99,CS}. Chavarriga and Garc\'{\i}a \cite{CG99}  obtained that
the origin is an isochronous center of the reversible system (1.4) with $\alpha_3=\alpha_4=0$ if and only if it can be transformed into  one of the following systems
\begin{eqnarray*}
\begin{aligned}
&CR_1:\ \dot{x}=y(-1+2\alpha_1x+2\alpha_2x^2), \ \ \dot{y}=x+\alpha_1(y^2-x^2)+2\alpha_2xy^2,\\
&CR_2:\ \dot{x}=y(-1+\alpha_1x+2\alpha_2x^2),  \ \dot{y}=x+\alpha_1y^2+2\alpha_2xy^2,\\
&CR_3:\ \dot{x}=y(-1+2\alpha_1x+3\alpha_2x^2-\alpha_2y^2),  \\
&\qquad \ \ \ \dot{y}=x+\alpha_1(y^2-x^2)-\alpha_2x^3+3\alpha_2xy^2,\\
&CR_4:\ \dot{x}=-y\big(-1+(\alpha_1-\beta_1)x+\beta_1(\alpha_1+\beta_1)(x^2-y^2)\big),\\
 &\qquad\ \ \ \dot{y}=x+\beta_1x^2+\alpha_1y^2+2\beta_1(\beta_1+\alpha_1)xy^2,\\
&CR_5:\ \dot{x}=-y(1-x)(1-2x),\ \ \dot{y}=x-2x^2+y^2+2x^3.
\end{aligned}
\end{eqnarray*}

The first integrals of the above systems are quite involved, and some of them are not even rational functions.  To our knowledge, there is  no  literature studying the limit cycle problems of these systems under polynomial small perturbations.

In the present paper, we intend to investigate the limit cycle problem of the cubic reversible  system $CR_5$ under perturbation of Li\'{e}nard type. That is,
\begin{eqnarray}
\dot{x}=-y(1-x)(1-2x),\ \ \dot{y}=x-2x^2+y^2+2x^3+\varepsilon f_1(x)y,
\end{eqnarray}
where $f_1(x)$ is a polynomial of degree $n$.
The change of variable $x=\frac{x_1+1}{2}$ allows us to rewrite the system (1.5) as
\begin{eqnarray}
 \dot{x}=xy-x^2 y,\ \ \dot{y}=y^2+\frac14x^3+\frac14x^2+\frac14x+\frac14+\varepsilon f(x)y,
\end{eqnarray}
where $0<|\varepsilon|\ll1$,$$f(x)=\sum\limits_{i=0}^na_ix^i,\ a_i\in\mathds{R}.$$
Here and below, we shall omit the subscript 1. System (1.6) for $\varepsilon=0$ has a first integral
\begin{eqnarray}
H(x,y)=\frac{1}{16}x^{-2}(x-1)^2(x^2+2x+1+4y^2)
\end{eqnarray}
with integrating factor $\mu(x)=\frac{1}{2}(x-1)x^{-3}$. The symmetric period annulus, surrounding the isochronous centre at $(x,y)=(-1,0)$ (corresponding to $h=0$), is denoted by $\Gamma_h=\{(x,y): H(x,y)=h,\, h>0\}$.

The total number of zeros of the following Abelian integral
\begin{eqnarray}I(h)=\frac12\sum\limits_{i=0}^na_i\oint_{\Gamma_h}(x-1)x^{i-3}ydx\end{eqnarray}
provides an upper bound for the number of limit cycles of system (1.6) bifurcating from the corresponding period annulus, and the existence of multiple simple zeros provides a lower bound of the number of limit cycles. Our main result is the following theorem.
\vskip 0.2 true cm

\noindent
{\bf Theorem 1.1}\, {\it  If the Abelian integral (1.8) is not identically zero, then, for any $n\geq1$, system (1.6) has at most $22n+6$ limit cycles bifurcating from the period annulus. Moreover, system (1.6) can have  one limit cycle when $n=1,2$ and two limit cycles when $n=3$ by choosing appropriate $a_i$.}
 \vskip 0.2 true cm

The present paper is built up as follows. We devote Section 2 to discuss the structure of Abelian integral which can be expressed by some elliptic integrals with rational functions, while Section 3 addresses the proof of our main result and numerical simulations are also presented.

\section{The algebraic structure of $I(h)$}
 \setcounter{equation}{0}
\renewcommand\theequation{2.\arabic{equation}}

The goal of this section is to obtain the algebraic structure of the Abelian integral $I(h)$. For abbreviation we use the notation
$$I_{i,j}(h)=\oint_{\Gamma_h}(x-1)x^{i-3}y^jdx,\ h\in(0,+\infty),$$
where $i,j\in\mathds{N}$. The orbit $\Gamma_h$ are symmetric with respect to the $x$-axis for $h\in(0,+\infty)$. Therefore, $I_{i,2j}(h)\equiv0$.

%In the sequel we will use the notation $P_k(h)$ to indicate the polynomial of $h$ with degree at most $k$.

\vskip 0.2 true cm

\noindent
{\bf Lemma 2.1}\, {\it For $k\in\mathds{N}$ and $k\geq2$, the following equality holds
\begin{eqnarray}\begin{aligned}
I_{k,3}(h)=&kI_{1,3}(h)+\lambda_1I_{0,3}(h)+\sum\limits_{i=2}^{k}P^i_1(h)I_{i,1}(h)\\&+\lambda_2I_{0,1}(h)+\lambda_3I_{1,1}(h)+\lambda_4I_{k+1,1}(h),
\end{aligned}\end{eqnarray}
where $P^i_1(h)$ are linear polynomials of $h$ for $i=2,\cdots,k$, and $\lambda_i\in\mathds{R}$ for $i=1,2,3,4$.}
 \vskip 0.2 true cm

\noindent
{\bf Proof}\, Applying  Green's formula, one obtains the relation
\begin{eqnarray}
\oint_{\Gamma_h}(x-1)^2x^iy^jdy=-\frac{i+2}{j+1}I_{i+3,j+1}(h)+\frac{i}{j+1}I_{i+2,j+1}(h).
\end{eqnarray}
 In order to establish the relations between the integrals $I_{i,j}(h)$, we differentiate $H(x,y)=h$ in (1.7) with respect to $x$ and multiply both sides by one-form $x^iy^{j-2}dx$. One finds that
\begin{eqnarray*}\begin{aligned}
I_{i,j}(h)=&-\frac14I_{i+3,j-2}(h)-\frac14I_{i+2,j-2}(h)-\frac14I_{i+1,j-2}(h)\\&-\frac14I_{i,j-2}(h)-\oint_{\Gamma_h}(x-1)^2x^{i-2}y^{j-1}dy.
\end{aligned}\end{eqnarray*}
In view of (2.2), one has the relation
\begin{eqnarray}\begin{aligned}
I_{i,j}(h)=&\frac{j}{i+j-2}\Big[\frac{i}{j}I_{i+1,j}(h)-\frac14I_{i+3,j-2}(h)-\frac14I_{i+2,j-2}(h)\\&-\frac14I_{i+1,j-2}(h)-\frac14I_{i,j-2}(h)\Big].
\end{aligned}\end{eqnarray}
Similarly, multiplying both sides of $H(x,y)=h$ by $(x-1)x^{i-3}y^jdx$ and integrating along $\Gamma_h$, we get another relation
\begin{eqnarray}\begin{aligned}
I_{i,j}(h)=&2(8h+1)I_{i-2,j}(h)-I_{i-4,j}(h)-4I_{i-2,j+2}(h)\\&+8I_{i-3,j+2}(h)-4I_{i-4,j+2}(h).
\end{aligned}\end{eqnarray}
Combining the equations (2.3) and (2.4), one easily obtains that
\begin{eqnarray}\begin{aligned}
I_{i,j}(h)=&\frac{i+j-3}{i+3j+1}\Big[\Big(16h+\frac{2i-10}{i+j-3}\Big)I_{i-2,j}(h)-I_{i-4,j}(h)-4I_{i-4,j+2}(h)\\&+
\frac{4i-4j-12}{i+j-3}I_{i-2,j+2}(h)-\frac{2j+4}{i+j-3}I_{i-1,j}(h)-\frac{2j+4}{i+j-3}I_{i-3,j}(h)\Big],
\end{aligned}\end{eqnarray}
and
\begin{eqnarray}\begin{aligned}
I_{i,j}(h)=&\frac{j}{i+j-1}\Big[\frac{1}{4}(16h+3)I_{i,j-2}(h)+\frac{i+3j-3}{j}I_{i-1,j}(h)\\
&+\frac{1}{4}I_{i+1,j-2}(h)+\frac{1}{4}I_{i-1,j-2}(h)-\frac{1}{4}I_{i-2,j-2}(h)-I_{i-2,j}(h)\Big].\end{aligned}\end{eqnarray}

Next we will prove (2.1) by induction on $k$. Indeed, a simple computation using (2.6) gives
\begin{eqnarray}\begin{aligned}
I_{2,3}(h)=&2I_{1,3}(h)-\frac34I_{0,3}(h)+\frac{3}{16}(16h+3)I_{2,1}(h)\\&
-\frac{3}{16}I_{0,1}(h)+\frac{3}{16}I_{1,1}(h)+\frac{3}{16}I_{3,1}(h),\\
I_{3,3}(h)=&3I_{1,3}(h)-\frac{27}{20}I_{0,3}(h)+\frac{3}{80}(64h+21)I_{3,1}(h)\\&
+\frac{3}{80}(144h+31)I_{2,1}(h)-\frac{27}{80}I_{0,1}(h)+\frac{3}{16}I_{1,1}(h)+\frac{3}{20}I_{4,1}(h),
\end{aligned}\end{eqnarray}
which imply that (2.1) holds for $k=2,3.$  Now assume that (2.1)  holds for $l\leq k-1$. Then, for $l=k$, one gets that
$$\begin{aligned}I_{k,3}(h)=&\frac{6}{k+2}\Big[\frac{k+6}{6}I_{k-1,3}(h)-\frac{1}{2}I_{k-2,3}(h)+\frac{1}{8}(1 6h+3)I_{k,1}(h)\\&
+\frac{1}{8}I_{k+1,1}(h)+\frac{1}{8}I_{k-1,1}(h)-\frac{1}{8}I_{k-2,1}(h)\Big],\end{aligned}$$
in view of (2.6).
This together with the induction hypothesis yields that (2.1) holds for $l=k.$ This completes the proof.\quad $\lozenge$
\vskip 0.2 true cm

 The recurrence formulas (2.5) and (2.6) play a key role in finding the detailed algebraic structure of the Abelian integral $I(h)$. The following lemma provides the necessary information in order to get the detailed expression of $I(h)$.
\vskip 0.2 true cm

\noindent
{\bf Lemma 2.2}\, {\it For $k\in\mathds{N}$ and $k\geq3$, there exist real polynomials $\bar{P}_{[\frac{k-3}{2}]}(h)$, $\bar{Q}_{[\frac{k-3}{2}]}(h)$, $\bar{R}_{[\frac{k-3}{2}]}(h)$ and $\bar{S}_{[\frac{k-2}{2}]}(h)$ of degrees at most $[\frac{k-3}{2}]$,  $[\frac{k-3}{2}]$, $[\frac{k-3}{2}]$ and $[\frac{k-2}{2}]$ such that
\begin{eqnarray}\begin{aligned}
&I_{k,1}(h)=\bar{P}_{[\frac{k-3}{2}]}(h)I_{0,1}(h)+\bar{Q}_{[\frac{k-3}{2}]}(h)I_{1,1}(h)+\bar{R}_{[\frac{k-3}{2}]}(h)I_{0,3}(h)+\bar{S}_{[\frac{k-2}{2}]}(h)I_{2,1}(h).
\end{aligned}\end{eqnarray}}

\noindent
{\bf Proof}\,  It follows from (2.6) that
\begin{eqnarray}\begin{aligned}
I_{1,3}(h)=&\frac{1}{4}(16h+3)I_{1,1}(h)+\frac{7}{3}I_{0,3}(h)+\frac{1}{4}I_{2,1}(h)\\&-I_{-1,3}(h)+\frac{1}{4}I_{0,1}(h)-\frac{1}{4}I_{-1,1}(h).
\end{aligned}\end{eqnarray}
A straightforward calculation using (2.5) gives that
\begin{eqnarray}\begin{aligned}
I_{3,1}(h)=&\frac{1}{7}(16h-4)I_{1,1}(h)-\frac{1}{7}I_{-1,1}(h)-\frac{4}{7}I_{-1,3}(h)\\
&-\frac{4}{7}I_{1,3}(h)-\frac{6}{7}I_{2,1}(h)-\frac{6}{7}I_{0,1}(h)\\
=&-I_{0,1}(h)-I_{1,1}(h)-\frac{4}{3}I_{0,3}(h)-I_{2,1}(h),
\end{aligned}\end{eqnarray}
on account of (2.9). It follows from (2.10) that (2.8) holds for $k=3$. Taking $(i,j)=(k,1)$ in (2.5), one obtains that
\begin{eqnarray}\begin{aligned}
I_{k,1}(h)=&\frac{1}{k+4}\big[4(k-4)I_{k-2,3}(h)-4(k-2)I_{k-4,3}(h)-(k-2)I_{k-4,1}(h)\\&+     (16(k-2)h+2k-10)I_{k-2,1}(h)-6I_{k-1,1}(h)-6I_{k-3,1}(h)\big]\\
\triangleq&\bar{P}_{[\frac{k-3}{2}]}(h)I_{0,1}(h)+\bar{Q}_{[\frac{k-3}{2}]}(h)I_{1,1}(h)+\bar{R}_{[\frac{k-3}{2}]}(h)I_{0,3}(h)+\bar{S}_{[\frac{k-2}{2}]}(h)I_{2,1}(h),
\end{aligned}\end{eqnarray}
using (2.1) and the induction hypothesis, where $\bar{P}_{[\frac{k-3}{2}]}(h)$, $\bar{Q}_{[\frac{k-3}{2}]}(h)$, $\bar{R}_{[\frac{k-3}{2}]}(h)$ and $\bar{S}_{[\frac{k-2}{2}]}(h)$ are polynomials of $h$ of degrees at most $[\frac{k-3}{2}]$,  $[\frac{k-3}{2}]$, $[\frac{k-3}{2}]$ and $[\frac{k-2}{2}]$. This completes the proof.\quad $\lozenge$\vskip 0.2 true cm

Inserting (2.8) into (1.11), and in light of (2.10), one immediately obtains the algebraic structure for the Abelian integral $I(h)$, as shown in the following lemma.

\vskip 0.2 true cm

\noindent
{\bf Lemma 2.3}\, {\it For $n\in\mathds{N}$ and $n\geq3$, there exist real polynomials ${P}_{[\frac{n-3}{2}]}(h)$, ${Q}_{[\frac{n-3}{2}]}(h)$, ${R}_{[\frac{n-3}{2}]}(h)$ and ${S}_{[\frac{n-2}{2}]}(h)$ of degrees at most $[\frac{n-3}{2}]$,  $[\frac{n-3}{2}]$, $[\frac{n-3}{2}]$ and $[\frac{n-2}{2}]$  such that
\begin{eqnarray}\begin{aligned}
&I(h)=P_{[\frac{n-3}{2}]}(h)I_{0,1}(h)+Q_{[\frac{n-3}{2}]}(h)I_{1,1}(h)+R_{[\frac{n-3}{2}]}(h)I_{3,1}(h)+S_{[\frac{n-2}{2}]}(h)I_{2,1}(h).
\end{aligned}\end{eqnarray}}

Let the complete elliptic integrals of first and second kind be
\begin{eqnarray}
E(k)=\int^{\frac \pi 2}_{0}\sqrt{1-k^2\sin^2(\theta)}d\theta,\ \ K(k)=\int^{\frac\pi 2}_0\frac{d\theta}{\sqrt{1-k^2\sin^2(\theta)}},\ k\in(-1,1).
\end{eqnarray}
Then, one can represent some generating integral functions of $I(h)$ as a combination of $E(k)$ and $K(k)$.
\vskip 0.2 true cm

\noindent
{\bf Lemma 2.4}\, {\it  The following equalities hold:
{\small\begin{eqnarray}\begin{aligned}
I_{1,1}(h)=&\frac{2}{\sqrt{8 h +1-4 \sqrt{4 h^{2}+h}}}E(k)-\frac{16h+2}{\sqrt{8h+1+4\sqrt{4h^2+h}}}K(k),\\
I_{3,1}(h)=&\frac{16h+2-8\sqrt{4h^2+h}\sqrt{8h+1-4\sqrt{4h^2+h}}}{3\sqrt{8h+1-4\sqrt{4h^2+h}}}K(k)\\&-\frac{16h+2}{3\sqrt{8h+1-4\sqrt{4h^2+h}}}E(k),\\
I_{0,1}(h)=&-4\pi h,\ \ I_{2,1}(h)=-4\pi h,
\end{aligned}\end{eqnarray}}where $${\small k=2 \sqrt{2}\, \sqrt{\frac{\sqrt{4 h^{2}+h}}{8 h +1+4 \sqrt{4 h^{2}+h}}},\ 0<k<1.}
$$}

\noindent
{\bf Proof}\, We begin with the proof by calculating the indefinite integral \begin{eqnarray}\int\sqrt{-x^4+(16h+2)x^2-1}dx.\end{eqnarray} For simplicity of the exposition, we introduce the notations
$$a=\frac{1}{\sqrt{8h+1-4\sqrt{4h^2+h}}},\ b=\frac{1}{\sqrt{8h+1+4\sqrt{4h^2+h}}},\ k=\frac{\sqrt{a^2-b^2}}{a}.$$
It is easy to check that $0<k<1$, in view of $h>0$.
The change of variable \begin{eqnarray}bx=\sqrt{1-k^2z^2},\end{eqnarray} allows us to rewrite the indefinite integral (2.15) as
\begin{eqnarray}
\frac{ak^4}{b^2}\int\frac{z^4-z^2}{\sqrt{(1-z^2)(1-k^2z^2)}}dz=\frac{ak^4}{b^2}(J_2-J_1),
\end{eqnarray}
where $$J_n=\int\frac{z^{2n}}{\sqrt{(1-z^2)(1-k^2z^2)}}dz,\ n=0,1,2.$$ In light of $$\frac{d\big(z\sqrt{(1-z^2)(1-k^2z^2)}\big)}{dz}=\frac{3k^2z^4-2(k^2+1)z^2+1}{\sqrt{(1-z^2)(1-k^2z^2)}},$$
one has that
\begin{eqnarray*}
J_2=\frac{1}{3k^2}\big(2(k^2+1)J_1-J_0+z\sqrt{(1-z^2)(1-k^2z^2)}\big).
\end{eqnarray*}
To summarize what we have proved that\begin{eqnarray*}\begin{aligned}
\int\sqrt{-x^4+(16h+2)x^2-1}dx=&\frac{2ak^2-ak^4}{3b^2}J_1-\frac{ak^2}{3b}J_0\\&+\frac{ak^2}{3b^2}z\sqrt{(1-z^2)(1-k^2z^2)}.
\end{aligned}\end{eqnarray*}
Thus, changing the variable of integration by letting $z=\sin\theta$, one gets that
\begin{eqnarray*}\begin{aligned}
I_{3,1}(h)=&\frac{2ak^2-ak^4}{3b^2}\int^1_0\frac{z^2}{\sqrt{(1-z^2)(1-k^2z^2)}}dz\\&-\frac{ak^2}{3b}\int^1_0\frac{1}{\sqrt{(1-z^2)(1-k^2z^2)}}dz\\
=&\frac{2a-ak^2-abk^2}{3b^2}K(k)-\frac{2a-ak^2}{3b^2}E(k),
\end{aligned}\end{eqnarray*}
and obtains the desired result on account of the relationship between $k$ and $h$.

In an analogous manner to that in (2.17), one has that
$$\begin{aligned}\int x^{-2}\sqrt{-x^4+(16h+2)x^2-1}dx=&-ak^2J_1-a(1-k^2)J_0\\&+a(1-k^2)\int\frac{1}{\sqrt{1-z^2}(1-k^2z^2)^\frac32}dz.\end{aligned}$$
Moreover, applying Newton-Leibniz's formula, one obtains that
$$\begin{aligned}I_{1,1}(h)=&-ak^2\int^1_0\frac{z^2}{\sqrt{(1-z^2)(1-k^2z^2)}}dz-a(1-k^2)\int^1_0\frac{1}{\sqrt{(1-z^2)(1-k^2z^2)}}dz\\&
+a(1-k^2)\int^1_0\frac{1}{\sqrt{1-z^2}(1-k^2z^2)^\frac32}dz.\end{aligned}$$
It can easily be checked that
$$\frac{d\Big(\frac{\sqrt{1-k^2z^2}\int^z_0\frac{\sqrt{1-k^2s^2}}{\sqrt{1-s^2}}ds-k^2z\sqrt{1-z^2}}{(1-k^2)\sqrt{1-k^2z^2}}\Big)}{dz}=\frac{1}{\sqrt{1-z^2}(1-k^2z^2)^\frac32}.$$
It follows from the above equality that
$$\int^1_0\frac{1}{\sqrt{1-z^2}(1-k^2z^2)^\frac32}dz=\frac{1}{1-k^2}\int^1_0\frac{\sqrt{1-k^2z^2}}{\sqrt{1-z^2}}ds.$$
A simple manipulation, using the change of variables $z=\sin\theta,$ leads to
$$I_{1,1}(h)=(ak^2-2a)K(k)+2aE(k).$$

We are now in a position to compute the remaining two integrals $I_{2,1}(h)$ and $I_{0,1}(h)$, which do not include the complete elliptic integrals of first and second kind.
Similar to calculating $I_{3,1}(h)$, using (2.16) we have that
{\small$$\begin{aligned}
&\int x^{-1}\sqrt{-x^4+(16h+2)x^2-1}dx=\frac{ak^4}{b}\int\frac{z^4-z^2}{\sqrt{1-z^2}(1-k^2z^2)}dz\\
=&-\frac{ak^2}{b}\int\frac{z^2}{\sqrt{1-z^2}}dz+\frac{a(k^2-1)}{b}\int\frac{1}{\sqrt{1-z^2}}dz-\frac{a(k^2-1)}{b}\int\frac{1}{\sqrt{1-z^2}(1-k^2z^2)}dz\\
=&\frac{ak^2-2a}{2b}\arcsin z+\frac{ak^2}{2b}z\sqrt{1-z^2}\\&+\frac{a\sqrt{1-k^2}}{b}\Big(\arctan\frac{1-kz-\sqrt{1-z^2}}{\sqrt{1-k^2}z}+\arctan\frac{1+kz-\sqrt{1-z^2}}{\sqrt{1-k^2}z}\Big)+C_1,
\end{aligned}$$}
where $C_1$ is an integration constant, which yields that
$$I_{2,1}(h)=\frac{ak^4}{b}\int^1_0\frac{z^4-z^2}{\sqrt{1-z^2}(1-k^2z^2)}dz=-4\pi h.$$

Some tedious computation leads to
{\small$$\begin{aligned}
&\int x^{-3}\sqrt{-x^4+(16h+2)x^2-1}dx=abk^4\int\frac{z^4-z^2}{\sqrt{1-z^2}(1-k^2z^2)^2}dz\\
=&ab\int\frac{1}{\sqrt{1-z^2}}dz+(2ab-abk^2)\int\frac{1}{\sqrt{1-z^2}(k^2z^2 -1)}dz\\&+ab(1-k^2)\int\frac{1}{\sqrt{1-z^2}(k^2z^2 -1)^2}dz\\
=&ab\arcsin z+\frac{abk^2z\sqrt{1-z^2}}{2(k^2z^2-1)}\\&-\frac{2ab-abk^2}{\sqrt{1-k^2}}\Big(\arctan\frac{1-kz-\sqrt{1-z^2}}{\sqrt{1-k^2}z}+\arctan\frac{1+kz-\sqrt{1-z^2}}{\sqrt{1-k^2}z}\Big)\\
&+\frac{ab(k^2-2)}{4\sqrt{1-k^2}}\Big(\arctan\frac{k-z}{\sqrt{(1-k^2)(1-z^2)}}-\arctan\frac{k+z}{\sqrt{(1-k^2)(1-z^2)}}\Big)+C_2,
\end{aligned}$$}
where $C_2$ is an integration constant, which leads to
$$I_{0,1}(h)=abk^4\int^1_0\frac{z^4-z^2}{\sqrt{1-z^2}(1-k^2z^2)^2}dz=-4\pi h.$$
The proof is completed. \quad $\lozenge$\vskip 0.2 true cm

The following proposition gives the detailed expression of the Abelian integral $I(h)$ which can be expressed by rational function and complete elliptic integrals with rational function coefficients.
\vskip 0.2 true cm

\noindent
{\bf Proposition 2.1}\, {\it When $n\geq3$, then the Abelian integral $I(h)$ can be written as
{\small\begin{eqnarray}
I(u)=\begin{cases}\frac{1}{u^{2[\frac{n-3}{2}]+3}}\Big[uP_{2[\frac{n-2}{2}]+2}(u^2)+Q_{4[\frac{n-3}{2}]+6}(u)K(k)+R_{2[\frac{n-3}{2}]+2}(u^2)E(k)\Big],\quad n\ odd,\\
\frac{1}{u^{2[\frac{n-2}{2}]+2}}\Big[P_{2[\frac{n-2}{2}]+2}(u^2)+Q_{4[\frac{n-3}{2}]+7}(u)K(k)+uR_{2[\frac{n-3}{2}]+2}(u^2)E(k)\Big],\quad n\ even,
\end{cases}\end{eqnarray}}
where the subscripts in polynomials $P$, $Q$ and $R$ represent their highest degrees, \begin{eqnarray}u=\sqrt{8h+1-4\sqrt{4h^2+h}},\ k=\sqrt{1-u^4},\ 0<u<1.\end{eqnarray}}

\noindent
{\bf Proof}\, Inserting (2.14) into (2.12) using (2.19) gives  (2.18).
This completes the proof.\quad $\lozenge$

\section{Proof of the main result}
 \setcounter{equation}{0}
\renewcommand\theequation{3.\arabic{equation}}

\subsection{Upper bound for the number of limit cycles}

In the sequel we will use the notation $\sharp\{u\in(\varsigma_1, \varsigma_2): \phi(u) = 0\}$ to indicate the number of zeros of the function $\phi(u)$ in the interval $(\varsigma_1, \varsigma_2)$, taking into account their multiplicities. In order to get the number of zeros of $I(h)$ for $h\in(0,+\infty)$, we need the following proposition, which can be proved by the Argument Principle based on the reference \cite{GLLZ}.
    \vskip 0.2 true cm

\noindent
{\bf Proposition 3.1}\, {\it Consider the function \begin{eqnarray}\Phi(z)=\varphi(z)K(z)+\psi(z)E(z),\end{eqnarray}
  where $\varphi(z)$ and $\psi(z)$ are polynomials of degrees at most $n$ and $m$, \begin{eqnarray}\begin{aligned}
&E(z)=\int^{\frac \pi 2}_{0}\sqrt{1-(1-z^4)\sin^2(\theta)}d\theta,\\ &K(z)=\int^{\frac\pi 2}_0\frac{d\theta}{\sqrt{1-(1-z^4)\sin^2(\theta)}},\ z\in(-1,1).
\end{aligned}\end{eqnarray} Then an upper bound for the number of zeros of $\Phi(z)$ for $z\in(-1,1)$, taking into account their multiplicities, is $n+m+1$.}
   \vskip 0.2 true cm

  \noindent
{\bf Proof}\, It is easy to check that $K(z)$ and $E(z)$ are analytic at $z=\pm1$ and $z=0$ is their singularity point. Then $K(z)$ and $E(z)$ can be continuously extended to single-valued analytic functions on the region
$$ D=\mathbb{C}\backslash \{z\in\mathbb{R},z>0\}.$$
Denote the upper and lower banks of the cut $\{z\in\mathbb{R},z>0\}$ by $L_1$ and $L_2$ respectively. Let $G_{R,\epsilon}\subset D$ be a simple connected region by removing a small disc $\{|z|<\epsilon\}$ and a real interval $[\epsilon,R]$ from $\{|z|\leq R\}$.

The following conclusions (i)-(iv) can be got from Lemmas 4 and 5 in \cite{GLLZ}, directly.

(i) The asymptotic expansions of $K(z)$ and $E(z)$ near 0 are
$$\begin{aligned}
&K(z)=2\ln2-2\ln z+O(z^4),\\ &E(z)=1-(\frac12-\ln4)z^4+z^4\ln z+O(z^8);
\end{aligned}$$

(ii) The asymptotic expansions of $K(z)$ and $E(z)$ near $+\infty$ are
$$\begin{aligned}
K(z)\sim z^{-2}\ln z,\ E(z)\sim z^2;
\end{aligned}$$

(iii) When $z\in L_i, i=1,2$,
$\textup{Im}\frac{E(z)}{K(z)}\neq0, \ \textup{Im}\frac{K(z)}{E(z)}\neq0;$
%$$K(u)E(u)\neq0,\ \ \textup{Im}\frac{K(u)}{E(u)}\neq0,\ \ \textup{Im}\frac{E(u)}{K(u)}\neq0,\ \ (\textup{Im}E(u))(\textup{Im}K(u))\neq0;$$

(iv) $K(z)$ and $E(z)$ have no zeros in the domain $D\subset \mathbb{C}$.

When $n\geq m+4$, by (iv), $$\sharp\{z\in D:I(z)=0\}=\{z\in D:\bar{I}(z)=\varphi(z)+\psi(z)\frac{E(z)}{K(z)}\}.$$
Our task now is to apply the Argument Principle to $\bar{I}(z)$ in the complex domain $G_{R,\epsilon}$ for $R$ and $1/\epsilon$ big enough.

By (ii), the number of complete turns along $\{|z|\leq R\}$ is at most $n$. By (iii), the number of complete turns of $\bar{I}(z)$ on $L_i, i=1,2$ is at most $m+1$, which includes less than one half turn on each bank. Finally, by (i), the number of complete turns of $\bar{I}(z)$ on $\{|z|=\epsilon\}$ goes to zero when $\epsilon\rightarrow0,$ owing to the fact that $\frac{E(z)}{K(z)}$ tends to a constant when $z\rightarrow0$. Hence, the total number of rotations of $\bar{I}(z)$ along the boundary of $G_{R,\epsilon}$ is at most $n+m+1$ and by the Argument Principle,
$$\sharp\{z\in D:\bar{I}(z)=0\}\leq n+m+1.$$

The case $n<m+4$ follows by the similar considerations by taking $\psi(z)+\varphi(z)\frac{K(z)}{E(z)}$  instead of $\varphi(z)+\psi(z)\frac{E(z)}{K(z)}$.
This completes the proof.\quad $\lozenge$
\vskip 0.2 true cm

\noindent
{\bf Lemma 3.1}\, {\it Let $k^2=1-z^4$ and
$$\Psi(z)=P_m(z)K(k)+Q_n(z)E(k),\ z\in(-1,1),$$
where $P_m(z)$ and $Q_n(z)$ are real polynomials of degrees $m$ and $n$. Then, the l$^{th}$-derivative of this expression is given by
$$\Psi^{(l)}(z)=\frac{1}{u^l(1-u^4)^l}\big(P_{m+4l}(z)K(k)+Q_{n+4l}(z)E(k)\big),$$
for $l\in\mathbb{N}$, where $P_{m+4l}(z)$ and $Q_{n+4l}(z)$ are real polynomials of degrees $m+4l$ and $n+4l$.}
 \vskip 0.2 true cm

  \noindent
{\bf Proof}\, The elliptic functions $K(k)$ and $E(k)$ satisfy the Picard-Fuchs equations
\begin{eqnarray*}
\frac{dK}{dk}=\frac{E-(1-k^2)K}{k(1-k^2)},\ \ \frac{dE}{dk}=\frac{E-K}{k},
\end{eqnarray*}
see \cite{BF}, formulas 710.00 and 710.02. The above two differential equations imply
\begin{eqnarray}
\frac{dK}{dz}=\frac{2z^4K-2E}{z(1-z^4)},\ \ \frac{dE}{dz}=\frac{2z^3(K-E)}{1-z^4},
\end{eqnarray}
in view of $k^2=1-z^4$. The conclusion follows directly by induction on $l$ using (3.3).
\vskip 0.2 true cm

  \noindent
{\bf Proof of Theorem 1}\, When $n\geq3$ is odd, in view of Proposition 2.1, we consider the number of zeros of the function
\begin{eqnarray}
\tilde{I}(u)=P_{4[\frac{n-2}{2}]+5}(u)+Q_{4[\frac{n-3}{2}]+6}(u)K(k)+R_{4[\frac{n-3}{2}]+4}(u)E(k)
\end{eqnarray}
instead of the zeros of $I(u)$, because $\tilde{I}(u)$ and $I(u)$ have the same number of zeros for $u\in(0,1)$. Thus, differentiating (3.4) $4[\frac{n-2}{2}]+6$ times using Lemma 3.1 gives that
\begin{eqnarray*}
\tilde{I}_1(u)=\frac{Q_{4[\frac{n-3}{2}]+16[\frac{n-2}{2}]+30}(u)K(k)+R_{4[\frac{n-3}{2}]+16[\frac{n-2}{2}]+28}(u)E(k)}{u^{4[\frac{n-2}{2}]+6}(1-u^4)^{4[\frac{n-2}{2}]+6}}.
\end{eqnarray*}
From Proposition 3.1 we obtain that the number of zeros of $\tilde{I}_1(u)$ in $(0,1)$ is bounded by $8[\frac{n-3}{2}]+32[\frac{n-2}{2}]+59$. Therefore, applying Rolle's Theorem $4[\frac{n-2}{2}]+6$ times, one has that
 $$\begin{aligned}\sharp\{u\in (0,1):{I}(u)=0\}&=\sharp\{u\in (0,1):\tilde{I}(u)=0\}\\&\leq 8[\frac{n-3}{2}]+36[\frac{n-2}{2}]+65\\&=22n-1.\end{aligned}$$

When $n\geq4$ is even, applying the method analogous to that used above, one can get that $$\begin{aligned}\sharp\{u\in (0,1):{I}(u)=0\}\leq 22n+6.\end{aligned}$$

When $n=1,2$, the proof of the result is quite similar to that given above and much easier and so is omitted.

\subsection{Lower bound for the number of limit cycles and numerical simulations}

In this section, we demonstrate that system (1.6) can have limit cycles for the given values of the coefficients of the perturbation polynomials as well as $|\varepsilon|$ small enough when $n=1,2,3$.
 Some numerical simulations are then presented showing that our results are right.

When $n=3$, we consider the following system
\begin{eqnarray}
 \dot{x}=xy-x^2 y,\ \ \dot{y}=y^2+\frac14x^3+\frac14x^2+\frac14x+\frac14+\sum\limits_{i=0}^3a_ix^iy,
\end{eqnarray}
where $$\begin{aligned}&a_0=\frac{1}{6\pi}\Big(8019E^2\big(\frac{4\sqrt{5}}{9}\big)-13020K\big(\frac{4\sqrt{5}}{9}\big)E\big(\frac{4\sqrt{5}}{9}\big) + 4121K^2\big(\frac{4\sqrt{5}}{9}\big)\Big)-1,\\
&a_1=-149K\big(\frac{4\sqrt{5}}{9}\big) - 219E\big(\frac{4\sqrt{5}}{9}\big),\\
&a_2=1,\ \ a_3=125K\big(\frac{4\sqrt{5}}{9}\big) - 405E\big(\frac{4\sqrt{5}}{9}\big).
\end{aligned}$$
Then computing the Abelian integral $I(u)$  for system (3.5) we have
\begin{eqnarray*}\begin{aligned}
I(u)=&-\frac{\pi(a_0 + a_2)}{8}(u^2 - 1)^2u^{-2} - \frac{a_1}{2}u^{-1}\big((u^4 + 1)K(\sqrt{1-u^4})-2E(\sqrt{1-u^4})\big) \\&
+ \frac{a_3}{6}u^{-3}\big((u^5 + u^4 - u + 1)K(\sqrt{1-u^4}) - (u^4 + 1)E(\sqrt{1-u^4})\big).
\end{aligned}\end{eqnarray*}
It is straightforward to show that the Abelian integral $I(u)$ has a double zero $u=\frac 13$ in the interval (0, 1) which corresponds a double zero $h=\frac{4}{9}$ of $I(h)$, see Fig. 2. It follows that for $|\varepsilon|\neq0$ small enough the  differential system (3.5) has a double limit cycle. Moreover we draw the  limit cycle that system (3.5) can exhibit when $\varepsilon= 10^{-8}$, see Fig. 3.

Taking $a_0=\frac{1}{6\pi}\Big(8019E^2\big(\frac{4\sqrt{5}}{9}\big)-13020K\big(\frac{4\sqrt{5}}{9}\big)E\big(\frac{4\sqrt{5}}{9}\big) + 4121K^2\big(\frac{4\sqrt{5}}{9}\big)\Big)$ and
the values of other parameters  are the same as those in system (3.5). Ay this moment, $I(u)$ has two simple zeros in (0,1), and two limit cycles exist in system (3.5), see Figs. 4 and 5.

When $n=2$, we consider system
\begin{eqnarray}
 \dot{x}=xy-x^2 y,\ \ \dot{y}=y^2+\frac14x^3+\frac14x^2+\frac14x+\frac14+\sum\limits_{i=0}^2a_ix^iy,
\end{eqnarray}
with $$\begin{aligned}a_0=-\frac{1}{2}E\big(\frac{\sqrt{15}}{4}\big),\ a_1=-\frac{9\pi}{128},\ a_2=\frac{17}{64}K\big(\frac{\sqrt{15}}{4}\big).
\end{aligned}$$
It can easily be verified that
\begin{eqnarray*}\begin{aligned}
I(u)=-\frac{\pi(a_0 + a_2)}{8}(u^2 - 1)^2u^{-2} - \frac{a_1}{2}u^{-1}\big((u^4 + 1)K(\sqrt{1-u^4})-2E(\sqrt{1-u^4})\big),
\end{aligned}\end{eqnarray*}
which has a zero $\frac12$ in (0,1), see Fig. 6, and a limit cycle of system (3.6) corresponds to it, see Fig. 7.

when $n=1$, a similar heuristic argument suggests that we consider the system
\begin{eqnarray}
 \dot{x}=xy-x^2 y,\ \ \dot{y}=y^2+\frac14x^3+\frac14x^2+\frac14x+\frac14+a_0y+a_1xy,
\end{eqnarray}
with $$\begin{aligned}a_0=25E\big(\frac{2\sqrt{6}}{5}\big)-13K\big(\frac{2\sqrt{6}}{5}\big),\ a_1=2\sqrt{5}\pi.
\end{aligned}$$
A straightforward computation yields that the Abelian integral $I(u)$ has a zero $\frac{\sqrt{5}}{5}$ which implies that system (3.7) has a limit cycle, see Figs. 8 and 9.

 \section{Conclusion}
 \setcounter{equation}{0}
\renewcommand\theequation{4.\arabic{equation}}

 The present study aims to obtain an upper bound  for the number of limit cycles for a cubic reversible system with an isochronous center  under the perturbations of $\varepsilon \sum\limits_{i=0}^na_ix^iy\frac{\partial}{\partial y}$. The algebraic structure of the Abelian integral is gotten by some iterative formulas of $I_{i,j}(h)$  and the bound is given by Argument Principle. Under our point of view, the main contribution of this article lies in the method to get the algebraic structure of the Abelian integral. Despite the previous proven result, there are still some open problems that should be solved for the cubic reversible system. However, it is difficult to study the limit cycle problem of $CR_5$ under general polynomial perturbations, because there are infinite generators for the corresponding Abelian integral, let alone $CR_i$ for $i=1,2,3,4.$ These difficulties may be overcome in the future by developing more powerful new approaches.
 \vskip 0.2 true cm

\noindent
{\bf CRediT authorship contribution statement}
 \vskip 0.2 true cm

\noindent
Jihua Yang: Conceptualization,  Funding acquisition, Investigation, Methodology, Resources, Supervision, Writing-original draft,	Writing-review $\&$ editing.
Qipeng Zhang: Validation, Writing-original draft.
\vskip 0.2 true cm

\noindent
{\bf Conflict of Interest}
 \vskip 0.2 true cm

\noindent
The authors declare that they have no known competing financial interests or personal relationships that could have
appeared to influence the work reported in this paper.
\vskip 0.2 true cm

\noindent
{\bf Data Availability Statement}
 \vskip 0.2 true cm

\noindent
No data was used for the research in this article. It is pure mathematics.
\vskip 0.2 true cm

\noindent
{\bf Acknowledgment}
 \vskip 0.2 true cm

\noindent
Supported by the National Natural Science Foundation of China(12161069) and the Ningxia Natural Science Foundation of China(2022AAC05044).

\end{document}